\documentstyle[11pt,leqno]{article}

\newtheorem{thm}{Theorem}[section]
\newtheorem{prop}[thm]{Proposition}
\newtheorem{lem}[thm]{Lemma}
\newtheorem{cor}[thm]{Corollary}

\def\ba{\begin{array}}
\def\ea{\end{array}}
\def\be{\begin{equation} \label}
\def\ee{\end{equation}}
\def\bea{\begin{eqnarray*}}
\def\eea{\end{eqnarray*}}
\def\beal{\begin{eqnarray} \label}
\def\eeal{\end{eqnarray}}
\def\bit{\begin{itemize}}
\def\eit{\end{itemize}}
\def\rf#1{(\ref{#1})}

\def\b{\beta}
\def\g{\gamma}
\def\Ga{\Gamma}
\def\d{\delta}
\def\D{\Delta}
\def\e{\varepsilon}
\def\L{\Lambda}

\def\s{\sigma}
 
\def\th{\vartheta}
\def\o{\omega}
\def\O{\Omega}

\def\LL{{\cal L}}

\def\F{{\cal F}}
\def\R{{\cal R}}
\def\T{{\cal T}}
\def\G{{\cal G}}

\def\Z{{\bf Z}}
\def\Gex{\mbox{\rm ex\,}\G}

\def\h{{\mbox{\tiny\rm hor}}}
\def\v{{\mbox{\tiny\rm vert}}}
\def\up{{\mbox{\tiny\rm up}}}
\def\lo{{\mbox{\tiny\rm down}}}
\def\le{{\mbox{\tiny\rm left}}}
\def\ri{{\mbox{\tiny\rm right}}}

\def\ls{\raisebox{0.4ex}{$\scriptstyle\leq$}}
\def\lss{\raisebox{0.4ex}{$\scriptstyle\leq$}*}

\begin{document}

\title{{\bf Percolation and number of phases\\ in the 2D Ising 
model}}
\author{ 
Hans-Otto Georgii\\ {\small\sl Mathematisches Institut der 
Universit\"at M\"unchen}\\ {\small\sl Theresienstr.\ 39, D-80333 
M\"unchen, Germany.} 
\and Yasunari Higuchi\\ {\small\sl Department of Mathematics, Faculty of 
Science}\\ {\small\sl Kobe University, Rokko, Kobe 657-8501, Japan.}
}

\date{}
\maketitle

\begin{quote}
We reconsider the percolation approach of Russo, Aizenman and Higuchi 
for showing that there exist only two phases in the Ising model on the 
square lattice. 
We give a fairly short alternative proof which is only based on FKG 
monotonicity and avoids the use of GKS-type inequalities 
originally needed for some background results. Our proof extends to 
the Ising model on other planar lattices such as the triangular 
and honeycomb lattice.
We can also treat the Ising antiferromagnet in a homogeneous field and 
the hard-core lattice gas model on $\Z^2$.
\end{quote}

\section{Introduction}

\medskip\noindent
One of the fundamental results on the two-dimensional ferromagnetic 
Ising model is the following theorem obtained independently 
in the late 1970s by Aizenman \cite{Aiz} and Higuchi \cite{Hig}
on the basis of the seminal work of Russo \cite{Rus}.

\medskip\noindent
{\bf Theorem}.
{\em For the ferromagnetic Ising model on $\Z^2$ with no external 
field and inverse temperature $\b>\b_c$, 
there exist only two distinct
extremal Gibbs measures $\mu^+$ and $\mu^-$.}

\medskip\noindent
The basic technique initiated by Russo  
consists of an interplay of three features of the Ising model: 
\bit
\item[--] the strong Markov property for random sets defined by 
geometric conditions involving clusters of constant spin, 

\item[--] the symmetry of the interaction under spin-flip and lattice 
automorphisms, and

\item[--] the ferromagnetic character of the interaction which 
manifests itself in FKG order and positive correlations. 
\eit
These ingredients led to a detailed understanding of the geometric 
features of typical configurations as described by the concepts of 
percolation theory.
In addition to these tools, the authors of \cite{Aiz,Hig,Rus} also 
needed the result that the
limiting Gibbs measure with $\pm$ boundary condition is a mixture of 
the two
pure phases. This result of Messager and Miracle--Sol\'e  \cite{MM} 
had been obtained by quite different means, namely some 
correlation inequalities of symmetry type in the spirit of GKS and 
Lebowitz inequalities. While 
such symmetry inequalities are a beautiful and powerful tool, they 
are quite different in character from the FKG inequality and have 
their own restrictions. It is therefore natural to ask whether 
Russo's random cluster method is flexible enough to prove the theorem 
without recourse to symmetry inequalities. On the one hand, this 
would allow to extend the theorem to models with less symmetries, 
while on the other hand one might gain a deeper understanding of 
possible geometric features of typical configurations.

In this paper we propose such a purely geometric reasoning which is 
only based on the three features above and avoids the use of the 
symmetry inequalities of Messager and Miracle-Sole \cite{MM}.
Despite this reduction of tools we could simplify the proof by 
an efficient combination of known geometric arguments. These include 
\bit 
\item [--] the Burton-Keane uniqueness theorem for infinite clusters 
\cite{BK},
\item [--] a version of Zhang's argument for the impossibility of 
simultaneous plus- and minus-percolation in $\Z^2$ (cf.\ Theorem 5.18 of 
\cite{GHM}),
\item [--] Russo's symmetry trick for simultaneous flipping of spins and 
reflection of the lattice \cite{Rus}, and 
\item [--] Aizenman's idea of looking at contour intersections in a 
duplicated system \cite{Aiz}.
\eit
We have tried to keep the paper reasonably 
self-contained, so that the reader will find a complete proof of the 
theorem.
As a payoff of the method we also obtain 
some generalizations. On the one hand, the arguments 
carry over to the Ising model on other planar lattices such as the 
triangular or the hexagonal lattice. On the other hand, in the case 
of the square lattice they cover 
also the antiferromagnetic Ising model in an external field as well 
as the hard-core lattice gas model.

\medskip\noindent
{\sl Acknowledgement. }We are indebted to Jeff Steif and an anonymous 
referee who asked for various additional details.

\section{Set-up and basic facts}
\label{facts}

Although we assume that the reader is familiar with the 
definition of the Ising model, let us start recalling a number of 
fundamental facts and introducing some notations. We assume 
throughout that the inverse temperature $\b$ exceeds the Onsager 
threshold $\b_{c}$, and that there is no external field, $h=0$.
The main ingredients we need are the following:
\smallskip

$\bullet$ 
the {\em configuration space\/} $\O=\{-1,1\}^{\Z^2}$, which is 
equipped 
with the Borel $\s$-algebra $\F$ and the local $\s$-algebras $\F_{\L}$
of events depending only on the spins in $\L\subset\Z^2$.

$\bullet$ 
the {\em Gibbs distributions} $\mu_{\L}^\o$ 
in finite regions $\L\subset\Z^2$ with boundary condition $\o\in\O$;
these enjoy the {\em Markov property} which says that 
$\mu_{\L}^\o(A)$ for
$A\in\F_{\L}$ depends only on the restriction of $\o$ 
to the boundary $\partial\L=\{x\not\in\L:|x-y|=1 
\mbox{ for some }y\in\L\}$ of $\L$, and the {\em finite energy 
property}, which states that $\mu_{\L}^\o(A)>0$ when $\emptyset\ne 
A\in\F_{\L}$. 

$\bullet$ 
the {\em Gibbs measures} $\mu$ on $(\O,\F)$ which, by 
definition,
satisfy $\mu(\cdot\,|\F_{\L^c})(\o) =\mu_{\L}^\o$ for $\mu$-almost 
all $\o$ and any finite $\L$;  we write $\G$ for the convex set 
of all Gibbs measures and $\Gex$
for the set of all extremal Gibbs measures.

$\bullet$ 
the {\em strong Markov property\/} of Gibbs measures,
stating that 
$\mu(\cdot\,|\F_{\Gamma^c})(\o) =\mu_{\Gamma(\o)}^\o$ for 
$\mu$-almost 
all $\o$ when $\Gamma$ is
any finite {\em random\/} subset of $\Z^2$ which is {\em determined from 
outside}, in that
$\{\Gamma=\L\}\in\F_{\L^c}$ for all finite $\L$, and
$\F_{\Gamma^c}$ is the $\s$-algebra of all events $A$ outside 
$\Gamma$, in the sense that 
$A\cap \{\Gamma=\L\}\in\F_{\L^c}$ for all finite $\L$. (Using the 
conventions $\mu_{\emptyset}^\o=\d_\o$ and $\F_{\emptyset^c}=\F$
we can in fact allow that $\Gamma$ takes the value $\emptyset$.)
For a proof  
one simply splits $\O$ into the disjoint sets $\{\Gamma=\L\}$ for 
finite $\L$.

$\bullet$ 
the {\em stochastic monotonicity} (or FKG order) of Gibbs 
distributions; writing $\mu\preceq\nu$ when 
$\mu(f)\leq\nu(f)$ for all increasing local (or, equivalently, 
all increasing bounded measurable) real functions $f$ on $\O$,
we have $\mu_{\L}^\o\preceq\mu_{\L}^{\o'}$ when $\o\leq\o'$, and  
$\mu_{\L}^\o\preceq\mu_{\D}^\o$ when $\D\subset\L$ and $\o\equiv 
+1$ on $\L\setminus\D$ (the opposite relation holds when $\o\equiv 
-1$ on $\L\setminus\D$).

$\bullet$ 
the {\em pure phases} 
$\mu^+,\,\mu^-\in\G$ obtained as limits for $\L\uparrow\Z^2$ of 
$\mu_{\L}^\o$ with $\o\equiv+1$ resp.\ $-1$, their invariance 
under all graph automorphisms of $\Z^2$,  
the sandwich relation 
$\mu^-\preceq\mu\preceq\mu^+$ for any other $\mu\in\G$, 
and the resulting extremality of $\mu^+$ and $\mu^-$.

$\bullet$ 
the characterization of extremal Gibbs measures by their 
{\em triviality on the tail $\s$-algebra} 
$\T=\bigcap\{\F_{\L^c}:\L\subset\Z^2\mbox{ finite}\}$; 
the fact that extremal Gibbs measures have {\em positive correlations}; 
and the {\em extremal 
decomposition} representing any Gibbs measure as the barycenter of a 
mass distribution on $\Gex$.

\smallskip\noindent
A general account of Gibbs measures can be found in \cite{Gii}, and 
\cite{GHM} contains an exposition of the 
Ising model and its properties related to stochastic monotonicity.

We will also use a class of transformations of $\O$ which preserve 
the Ising Hamiltonian, and thereby the class $\G$ of Gibbs measures. 
These transformations are
\smallskip

$\bullet$
the {\em spin-flip transformation} 
$T:\o=(\o(x))_{x\in\Z^2}\to (-\o(x))_{x\in\Z^2}\,;$

$\bullet$ the {\em translations} $\th_{x}$, $x\in\Z^2$, which are 
defined by $\th_{x}\o(y)=\o(y-x)$ for $y\in \Z^2$, and in particular 
the horizontal and vertical shifts 
$\th_{\h}=\th_{(1,0)}$ resp.\  $\th_{\v}=\th_{(0,1)}$; and

$\bullet$ the {\em reflections\/} in lines $\ell$ through lattice 
sites: for any $k\in\Z$ we write
\[
R_{k,\h}: \Z^2\ni x=(x_1,x_2) \to (x_1,2k-x_2)\in\Z^2
\]
for the reflection in the horizontal line $\{x_2=k\}$,
and similarly $R_{k,\v}$ for the reflection in the vertical line
$\{x_1=k\}$. For $k=0$ we simply write $R_\h=R_{0,\h}$ and 
$R_\v=R_{0,\v}$. All these reflections act canonically on $\O$.
\medskip

We will investigate the geometric behavior of typical configurations 
in {\em half-planes\/} of $\Z^2$. These are sets of the form 
\[
\pi=\{ x=(x_1, x_2) \in \Z^2 :\ x_i \ge k\}
\]
with $k\in\Z$, $i\in\{1,2\}$, or with `$\ge$' replaced by `$\leq$'.
The line $\ell =\{ x\in \Z^2: \ x_i= k \}$ is called the associated 
{\em boundary line}. In particular, we will consider 
\bit
\item the upper half-plane 
$ \pi_\up = \{ x=(x_1, x_2) \in \Z^2: \ x_2 \geq 0 \} $,
\item the downwards half-plane
$ \pi_\lo = \{ x = (x_1, x_2) \in \Z^2: \ x_2 \leq 0 \} $,
\eit
and the analogously defined right half-plane $ \pi_\ri$
and left half-plane $\pi_\le$. We will also work with
\bit
\item the left horizontal semiaxis
$\ell_\le = \{x=(x_1,x_2) \in \Z^2: \ x_1 \leq 0,\, x_2=0\}$, and
\item the right semiaxis
$ \ell_\ri = \{x = (x_1, x_2) \in \Z^2: \ x_1 \geq 0,\, x_2 = 0 \}$.
\eit
\medskip

In the rest of this section we state three fundamental results on 
percolation in the Ising model. By the symmetry between 
the spin values $+1$ and $-1$, 
these results also hold when `$-$' and `$+$' are 
interchanged. Similarly, all notations introduced with one 
sign will be used accordingly for the opposite sign. 

We first recall some basic concepts of percolation theory.
A finite {\em path} is a sequence  $p=(x_1, x_2,\ldots,x_k)$ 
of pairwise distinct lattice points such that, for any 
$i\in \{2,\ldots, k\}$, $x_{i-1}$ and $x_i$ are nearest 
neighbors (i.e., have Euclidean distance 1). The number $k$ is called
the {\em length\/} of $p$, and $x_1$ and $x_k$ are its starting 
resp.\ final point. 
A path $p$ is called a {\em path in a subset\/} 
$S\subset\Z^2$ if all $x_i$ belong to $S$. We say that $p$ {\em 
meets\/} or {\em touches\/} $S$ if some $x_i$ is contained in $S$ or 
a nearest neighbor of a point in $S$. We will also speak of
infinite paths $(x_1, x_2, \ldots)$ and doubly infinite 
paths $(\ldots, x_{-1},x_0, x_1, \ldots)$ in the obvious sense. 
A path $p$ is called a {\em circuit\/} if $x_1$ and $x_k$ 
are nearest neighbors, and a {\em semicircuit\/} in a half-plane $\pi$ 
if it is contained in $\pi$ and $x_1$ and $x_k$ belong to the boundary 
line of $\pi$. A region $C\subset\Z^2$ is called {\em connected\/} 
if for any $x,y\in C$ there exists a path in $C$
from $x$ to $y$. A {\em cluster\/} in a region $S\subset\Z^2$ is a 
maximal connected subset $C$ of $S$. It is called infinite if its 
cardinality is infinite. Infinite clusters will be denoted by the 
letter $I$, with suitable sub- and superscripts. 

Given any configuration $\o\in\O$, we consider the set 
$S^+(\o)=\{x\in\Z^2: \o(x)=+1\}$ of $+$ spins. A path
(resp.\ circuit, semicircuit, cluster) in $S^+(\o)$ is called a 
{\em $+$path\/} (resp.\ {\em $+$circuit, $+$semicircuit, $+$cluster\/}) 
for $\o$, and two points $x,\,y$ are said to be {\em $+$connected\/} if there 
exists a $+$path from $x$ to $y$. 

We also need to work with the conjugate graph structure on $\Z^2$, 
for which two points are considered as neighbors if their Euclidean 
distance is either 1 or $\sqrt 2$, i.e., if they are either nearest 
neighbors or diagonal neighbors. This graph structure is indicated by 
a star and leads to the concepts of $*$paths, $*$circuits, 
$*$semicircuits, $*$connectedness, $*$clusters, $+*$paths, 
$+*$semicircuits, and so on. Note that each path is a fortiori a 
$*$path, and each cluster is contained in some $*$cluster.

The starting point of the random cluster method is the following result 
of \cite{CNPR,Rus}. Let $E^+$ denote the 
event that there exists an infinite $+$cluster $I^+$ in 
$\Z^2$, and define $E^-$, $E^{+*}$, $E^{-*}$ analogously.
Note that $E^+\subset E^{+*}$ and $E^-\subset E^{-*}$.
(Throughout this paper we will use the letter $E$ to denote events 
concerning existence of infinite clusters.) 
\begin{lem}{\bf (Existence of infinite clusters)}\label{lem:CNPR}
If $\mu\in\G$ is different from $\mu^-$, there exists with 
positive probability an infinite $+$cluster.
That is, $\mu(E^+)>0$ when $\mu\ne\mu^-$.
\end{lem}
{\sl Proof}. 
Suppose that $\mu(E^+)=0$. Then any given square 
$\D$ is almost surely surrounded by a $-*$circuit, and with 
probability close to $1$ such a circuit can already be found within a 
square $\L\supset\D$ provided $\L$ is large enough. If this
occurs, we let $\Gamma$ be 
the largest random subset of $\L$ which is the 
interior of such a $-*$circuit. (A largest such set exists because 
the union of such sets is again the interior of a $-*$circuit.)
In the alternative case we set $\Gamma=\emptyset$. By 
maximality, $\Gamma$ is determined from outside. The strong Markov 
property together with the stochastic monotonicity 
$\mu^-_\Gamma\preceq \mu^-$ therefore implies (in the limit $\L\uparrow\Z^2$) 
that $\mu\preceq\mu^-$ on $\F_\D$.
Since $\D$ was arbitrary and $\mu^-$ is minimal we find that 
$\mu=\mu^-$, and the lemma is proved. $\Box$

\medskip\noindent
The next lemma is a variant of another result of Russo \cite{Rus}. 
\begin{lem}{\bf (Flip-reflection domination)}\label{TR}
Let $\mu \in \G$ and $R$ any reflection, and suppose that for
$\mu$-almost all $\o$ each finite $\D\subset \Z^2$ is surrounded by 
an $R$-invariant 
$*$circuit $c$ such that $\o \ge R\circ T(\o)$ on $c$.
Then $\mu \succeq \mu \circ R \circ T$.
\end{lem}
{\sl Proof}. Another way of stating the assumption 
is that for
any finite $R$-invariant $\D$ and $\mu$-almost all $\o$ 
there exists a finite $R$-invariant random set $\Ga(\o)\supset\D$ 
such that $\o \ge R\circ T(\o) \mbox{ on }\partial\Ga(\o)$.
Given any $\e>0$, we can thus find an $R$-invariant $\L$ so large that
with probability at least $1-\e$ such 
an $R$-invariant $\Ga(\o)$ exists within $\L$.
Since the union of any two such $\Ga(\o)$'s enjoys the same 
properties,
we can assume that $\Ga(\o)$ is chosen maximal in $\L$; in the case 
when no 
such $\Ga(\o)$ exists we set $\Ga(\o)=\emptyset$. 
The maximality of $\Ga$ implies that the events $\{\Ga = G\}$ are 
measurable with respect to $\F_{\L\setminus G}$.
For any increasing ${\cal F}_{\D}$-measurable  function $f\geq 0$
we thus get from the strong Markov property
$$
\mu(f) \geq \mu\bigg( \mu_\Ga^\cdot( f)\; 1_{\{\Ga \ne\emptyset\}}\bigg)\;.
$$
However, if $\Ga(\o)\ne\emptyset$ then
$\o \ge R\circ T(\o) \mbox{ on }\partial \Ga(\o)$. By stochastic 
monotonicity, for such $\o$ we have 
$$
\mu_{\Ga(\o)}^\o( f )\ge 
\mu_{\Ga(\o)}^{R\circ T(\o)}( f )=\mu_{\Ga(\o)}^\o( f\circ R\circ T 
)\;,
$$
where the identity follows from the 
$R$-invariance of $\Ga$ and the $R\circ T$ -invariance of the 
interaction. Hence
$$
\mu(f)\geq \mu\bigg( f\circ R\circ T \;1_{\{\Ga \ne \emptyset\}}\bigg)
\geq \mu ( f\circ R\circ T ) - \e\,\| f \|_\infty \; .
$$
The lemma thus follows by letting $ \e \to 0 $ and $\D\uparrow\Z^2$.
 $\Box$

\medskip\noindent
A third useful result of Russo \cite{Rus} is the following. To state 
it we need to introduce two notations. First, 
let 
\[
\theta=\mu^+(0 \in I^{+*})
\]
be the $\mu^+$-probability that the origin belongs to 
an infinite $+*$cluster. Lemma \ref{lem:CNPR} implies that $\theta>0$.
Secondly, for a half-plane $\pi$ with boundary line $\ell$ and a 
$*$semicircuit $\s$ in $\pi$ we write $\mbox{Int\,}\s$ for the unique 
subset of $\Z^2$ which is invariant under the reflection $R$ in 
$\ell$ 
and satisfies $\pi\cap\partial(\mbox{Int\,}\s)=\s$; we call 
$\mbox{Int\,}\s$ the interior of $\s$. 
\begin{lem} {\bf(Point-to-semicircuit lemma)}
\label{point-to-semicircuit}
Let $\pi$ be some half-plane with boundary line $\ell$, $x\in\ell$, 
and $\s$ a $*$semicircuit in $\pi$ with interior 
$\L=\mbox{\rm Int\,}\s\ni x$. Let $\o\in\O$ be such that
$\o\equiv +1$ on $\s$. Then 
\[
\mu_\L^\o\bigg( x \mbox{ \rm is in $\L$ $+*$connected to $\s$} \bigg) \geq 
\theta/2\;.
\]
\end{lem}
{\sl Proof}. By stochastic monotonicity we can assume that $\o\equiv 
-1$ on $\partial\L\setminus\s$. We then have $\o\geq R\circ T(\o)$ on 
$\partial\L$, and therefore
$\mu_\L^\o\succeq\mu_\L^\o\circ R\circ T$. To exploit this relation
we let $B_{x,\s}$ be the event that there exists a $+*$path in $\L$ 
from $x$ to $\s$, $C_{x,\s}$ the event that $x$ is surrounded by a 
$+*$circuit in $\L\cup\s$ which is $+*$connected to $\s$, and 
$D_{x,\s}=B_{x,\s}\cup C_{x,\s}$. 
Then 
\be{lasso}\mu_\L^\o(D_{x,\s}\cup R\circ T(D_{x,\s}))=1\:. 
\ee
Indeed, suppose
that $\o(x)=+1$, but $B_{x,\s}$ does not occur. Then the $+*$cluster
containing $x$ does not meet $\s$. Its outer boundary belongs to a
$-*$cluster, which either touches $R(\s)$ so that $R\circ T(C_{x,\s})$
occurs, or not --- in which case we consider the $+*$cluster containing
its outer boundary, and so on. After finitely many steps we see that either
$C_{x,\s}$ or $R\circ T(C_{x,\s})$ must occur. \rf{lasso} is an immediate consequence. It follows
that $\mu_\L^\o(D_{x,\s})\geq1/2$. Hence 
$$\mu_\L^\o(B_{x,\s}) 
\ge\mu_\L^\o(B_{x,\s}|C_{x,\s})\,\mu(D_{x,\s})
\ge \mu_\L^\o(B_{x,\s}|C_{x,\s})/2\,.
$$
But if $C_{x,\s}$ occurs then there exists a largest random set 
$\Ga\subset\L$ containing $x$ such that $\partial\Ga$ forms a 
$+*$circuit and is $+*$connected to $\s$. Writing $B_{x,\partial\Ga}$ 
for the event that $x$ is $+*$connected to $\partial\Ga$ and
using the strong Markov property we thus find that
$$
\mu_\L^\o(B_{x,\s}|C_{x,\s})= 
\mu_\L^\o\,(\,\mu^+_{\Gamma}(B_{x,\partial\Ga})\,|\,C_{x,\s}\,)
\geq\theta
$$
because $\mu^+_{\Gamma}(B_{x,\partial\Ga})\geq\theta$ by 
stochastic monotonicity. Together with the previous inequality this 
gives the result. $\Box$

\section{Percolation in half-planes} 

In this section we will prove that there exist plenty of infinite 
clusters of constant spin in the half-planes of $\Z^2$. 
In particular, this 
will show that all translation invariant $\mu\in{\G}$ 
are mixtures of $\mu^+$ and $\mu^-$.
We will use two pearls of percolation theory, 
the Burton-Keane uniqueness theorem for infinite clusters \cite{BK}, 
and Zhang's argument for 
the non-existence of two infinite clusters of opposite sign in $\Z^2$.
(In the present context, these two results were obtained first
in \cite{GKR}.)

For a given half-plane $\pi$ we let $E^+_\pi$ denote the event that 
there exists an infinite $+$cluster in $\pi$. When this occurs, we 
will write $I^+_\pi$ for such an infinite $+$cluster in $\pi$. (As we 
will see, such clusters are unique, so that this notation does not 
lead into conflicts.) In case of the standard half-planes,
we will only keep the 
directional index and omit the $\pi$; for example, we write
$E^+_{\up}$ for $E^+_{\pi_{\mbox{\tiny up}}}$.
Similar notations will be used for $+*$clusters 
and for the sign $-$ instead of $+$.

Let us say that $(\pi,\,\pi')$ is a pair of {\em conjugate 
half-planes\/} if $\pi$ and $\pi'$ share only a common boundary line. 
An associated
pair $(I^+_{\pi},I^+_{\pi'})$ or $(I^-_{\pi},I^-_{\pi'})$ of infinite 
clusters of the same sign in $\pi$ and $\pi'$ will be called an
{\em infinite butterfly}. (This name alludes to the assumption that 
the two infinite `wings' have the same `color', but is not meant to 
suggest that they are symmetric and connected to each other, 
although the latter will turn out to be true.)  We will say that a 
statement holds $\G$-almost surely if it holds $\mu$-almost surely 
for all $\mu\in\G$.
\begin{lem}{\bf(Butterfly lemma)}\label{butterfly}
$\G$-almost surely there exists at least one infinite butterfly.
\end{lem}
{\bf Proof: } Suppose the contrary. By the extremal decomposition 
theorem
and the fact that the existence of infinite butterflies is a tail 
measurable event, there is then some $\mu\in{\Gex}$ for which there
exists no infinite butterfly $\mu$-almost surely.
We will show that this is impossible.

{\em Step 1. }First we observe that $\mu$ is $R\circ T$-invariant 
for all reflections $R=R_{k,\h}$ or $R_{k,\v}$, 
and in particular is periodic under translations.
Indeed, let $(\pi,\,\pi')$ be conjugate half-planes with common 
boundary line $\ell$ and $R$ the 
reflection in $\ell$ mapping $\pi$ onto $\pi'$. By the absence of 
infinite butterflies, at least
one of the half-planes $\pi$ and $\pi'$ contains no infinite
$-$cluster, and this or the other half-plane contains no infinite
$+$cluster. In view of the tail triviality of $\mu$, we can assume 
that $\mu(E^-_{\pi})=0$. This means that for $\mu$-almost all $\o$ 
every finite $\D\subset\pi$ is surrounded 
by some $+*$semicircuit $\g$ in $\pi$. For such a $\g$, $c=\g\cup 
R(\g)$
is an $R$-invariant $*$circuit that surrounds $\D\cup R(\D)$
and satisfies $\o \ge R\circ T(\o)$ on $c$. By Lemma \ref{TR}, 
this gives the flip-reflection domination
$\mu \succeq \mu \circ R \circ T$. Since also $\mu(E^+_{\pi})=0$ or
$\mu(E^+_{\pi'})=0$, we conclude in the same way that
$\mu \preceq \mu \circ R \circ T$, so that $\mu=\mu\circ R\circ T$. 
Since both $\th_\h^2$ and $\th_\v^2$ are compositions of two 
reflections, the invariance under the translation group 
$(\th_x)_{x\in 2\Z^2}$ follows.

{\em Step 2. }We now take advantage of
the Burton--Keane uniqueness theorem \cite{BK}, stating that for 
every periodic $\mu$ with finite energy there exists at most
one infinite $+$ (resp.\ $-$) cluster, and Zhang's symmetry argument
(cf.\ \cite{GHM}, Theorem 5.18) 
deducing from this uniqueness the absence of simultaneous $+$ and 
$-$percolation. (In reference \cite{BK}, the uniqueness of the 
infinite cluster is only stated for translation invariant $\mu$, but 
the argument works in the same way by applying the ergodic theorem 
to the subgroup $(\th_x)_{x\in 2\Z^2}$. It is also not shown 
there that the finite energy property remains valid under 
ergodic decomposition. Although this follows from 
Theorem (14.17) of \cite{Gii} in the present setting, and a similar 
argument in general, we do not need this here because our $\mu$ is 
extremal, and therefore $(\th_x)_{x\in 2\Z^2}$-ergodic by 
Proposition (14.9) of \cite{Gii}.)

We start noting that, by the flip-reflection symmetry of $\mu$, 
$\mu$ is different
from $\mu^+$ and $\mu^-$, so that by Lemma \ref{lem:CNPR}, the 
tail triviality of $\mu$, and
the Burton--Keane uniqueness theorem there exist both a unique 
infinite 
$+$cluster $I^+$ and a unique infinite $-$cluster $I^-$ 
in the whole plane $\Z^2$ $\mu$-almost surely. 
We now choose a square $\L=[-n,n]^2\cap\Z^2$ so large that
$\mu(\L\cap I^+\ne\emptyset)>1-2^{-12}$.
Let $\partial_k\L$ be the intersection of
$\partial\L$ with the $k$'th quadrant, and let $A^+_{k}$ 
be the increasing event that there exists an 
infinite $+$path in $\L^c$ starting from some site 
in $\partial_k\L$. Define $A^-_{k}$ analogously. Since 
\[
\{\L\cap I^+\ne\emptyset\}\subset \bigcup_{k=1}^4 A^+_{k}
\]
and $\mu$ (as an extremal Gibbs measure) has positive correlations,
it follows that
\[
\prod_{k=1}^4 \mu(\O\setminus A^+_{k})\leq \mu(\bigcap_{k=1}^4
\O\setminus{A^+_{k}})
\leq\mu(\L\cap I^+=\emptyset) < 2^{-12}\;,
\]
whence there exists some $k\in\{1,\ldots,4\}$ such that 
$\mu(\O\setminus{A^+_{k}})<2^{-3}$. For notational convenience 
we assume that $k=1$. 
By the flip-reflection symmetry shown above, we find that
\[
\mu(A^+_{1}\cap A^-_{2} \cap A^+_{3} \cap A^-_{4})>1-4\cdot 
2^{-3}=1/2\;,
\]
which is impossible because on this intersection the infinite clusters
$I^+$ and $I^-$ cannot be both unique. This contradiction concludes 
the proof of the lemma. $\Box$

\medskip\noindent
The butterfly lemma leads immediately to the following result first 
obtained by Messager and Miracle-Sole \cite{MM} by means of 
correlation inequalities of symmetry type; the following proof 
appeared first in \cite{GHM}.
\begin{cor}\label{cor:MM}
{\bf (Periodic Gibbs measures) }
Any periodic $\mu\in\G$ is a mixture of $\mu^+$ and $\mu^-$.
\end{cor}
{\sl Proof}. Suppose $\mu\in\G$ is invariant under 
$(\th_x)_{x\in p\Z^2}$ for some period $p\geq 1$. Conditioning $\mu$
on any periodic tail event $E$ we obtain again a periodic Gibbs measure.
It is therefore sufficient to show that $\mu(E^+\cap E^-)=0$. Indeed,
the butterfly lemma then shows that $\mu(E^+) +\mu(E^-)=1$, and Lemma 
\ref{lem:CNPR} implies that $\mu(\,\cdot\,|E^+)=\mu^+$ and 
$\mu(\,\cdot\,|E^-)=\mu^-$ whenever these conditional probabilities 
are defined. Hence $\mu=\mu(E^+)\,\mu^+ +\mu(E^-)\,\mu^-$.

Suppose by contraposition that $\mu(E^+\cap E^-)>0$. Since 
$E^+\cap E^-$ is invariant and tail measurable, we can in fact 
assume that $\mu(E^+\cap E^-)=1$; otherwise we replace $\mu$ by 
$\mu(\,\cdot\,|E^+\cap E^-)$. By the butterfly lemma,
there exists a pair $(\pi,\pi')$ of conjugate halfplanes, say 
$\pi_\up$ and $\pi_\lo$, and a sign, say $+$, such that both 
half-planes contain an infinite clusters of this sign with 
positive probability. Since $\mu(E^-)=1$ by assumption, 
we can find a large square $\D$ such that with 
positive probability $\D$ meets infinite $+$clusters in $\pi_\up$ and 
 $\pi_\lo$ and also an infinite $-$cluster. This $-$cluster leaves 
$\D$ either on the left or on the right between the two infinite 
$+$clusters. We can assume that the latter occurs with positive 
probability. By the
finite energy property, it then follows 
that also $\mu(A_0)>0$, where for $k\in p\Z$ we write
$A_k$ for the event that the point $(k,0)$ belongs to a two-sided 
infinite $+$path with its two halves staying in $\pi_\up$ resp.\ 
$\pi_\lo$, and $(k+1,0)$ belongs to an infinite $-$cluster.

Let $A$ be the event that $A_k$ occurs for infinitely many $k<0$ and
infinitely many $k>0$. The horizontal periodicity and Poincar\'e's
recurrence theorem (cf.\ Lemma (18.15) of \cite{Gii}) then show that
$\mu(A_0\setminus A)=0$, and therefore $\mu(A)>0$. But on $A$ there 
exist infinitely many $-$clusters which are separated from each other 
by the infinitely many `vertical' $+$paths. This contradicts the 
Burton--Keane theorem.
$\Box$

\medskip\noindent
The preceding argument actually shows that $\mu(E^{-*}\cap 
E^{+*})=0$ whenever $\mu\in\G$ is periodic. Since $\mu^+(E^+)=1$ by 
Lemma \ref{lem:CNPR} and tail triviality, this shows that in the 
$+$phase the $+$spins 
form an infinite sea with only finite islands.
\begin{cor}\label{ocean}
{\bf (Plus-sea in the plus-phase)}
$\mu^+(E^{-*})=0$. Hence, 
 $\mu^+$-almost surely there exists a unique infinite 
$+$cluster $I^+$ in $\Z^2$ which surrounds each finite set. 
\end{cor}
We note that in contrast to Zhang's argument 
(cf.\ Theorem 5.18 of \cite{GHM})
our proof of the preceding corollary does not rely on 
the reflection invariance of $\mu^+$ but only on its periodicity, and 
thus can be extended to the setting of Section \ref{extensions} below.

We conclude this section with the observation that percolation in 
half-planes is not affected by spatial shifts.
\begin{lem}\label{shift-lemma}
{\bf (Shift lemma)}
Let $\pi$ and $\tilde\pi$ be two half-planes such that 
$\pi\supset\tilde\pi$,
i.e., $\pi$ and $\tilde\pi$ are translates of each other. Then
$E^+_{\pi}= E^+_{\tilde\pi}$ $\G$-almost surely, and similarly with 
$-$ instead of $+$.
\end{lem}
{\sl Proof}. Since trivially $E^+_{\pi}\supset E^+_{\tilde\pi}$,
we only need to show that $E^+_{\pi}\subset E^+_{\tilde\pi}$ 
$\G$-almost surely. For definiteness we consider the case when 
$\pi=\pi_\up=\{x_2\ge0 \}$ and 
$\tilde\pi =\{x_2\ge 1\}$. Take any $\mu\in\Gex$, and 
suppose that $\mu( E^+_{\tilde\pi})=0$. Then for almost all $\o$ and 
any $n\geq1$ 
there exists a smallest $-*$semicircuit $\s_n(\o)$ in $\tilde\pi$ 
containing $\D_n\cup\s_{n-1}(\o)$ in its interior; 
here $\D_n= [-n,n]\times[1,n]$ and $\s_0=\emptyset$. 
Let $x_n(\o)\in\ell_\le$ and $y_n(\o)\in\ell_\ri$ be the two points
facing the two endpoints 
of $\s_n(\o)$; these are $\F_{\tilde\pi}$-measurable functions of 
$\o$, and the random sets $\{x_n,y_n\}$ are pairwise disjoint. Let 
$A_n$ 
be the event that the spins at  $x_n$ and $y_n$ take value $-1$. 

We claim that $A_n$ occurs for infinitely many $n$ with probability 
1. Indeed, fix any $N\ge1$, $x\in\ell_\le$, $y\in\ell_\ri$ and 
let $B_{N,x,y}=\{x_N=x,y_N=y\}\cap\bigcap_{n>N}A_n^c$. Then we can 
write
\[
\mu(A_{N}\cap B_{N,x,y})
=\mu\bigg(\mu_{\{x,y\}}^\cdot(\o(x)=\o(y)=-1)\,1_{B_{N,x,y}}\bigg)
\ge \d^2\, \mu(B_{N,x,y})
\]
because $B_{N,x,y}$ only depends on the configuration outside 
$\{x,y\}$,
and the one-point conditional probabilities of $\mu$
are bounded from below by $\d=[ 1 + e^{8\beta } ]^{-1}$.
Summing over $x,y$ we obtain
$\mu(\bigcap_{n\geq N}A_n^c)\leq (1-\d^2)\,\mu(\bigcap_{n> N}A_n^c)$, 
and iteration gives
$\mu(\bigcap_{n\ge N}A_n^c)=0$. Letting $N\to\infty$ we get the claim.

We now can conclude that 
with probability $1$ each box $[-n,n]\times[0,n]$ is surrounded
by a $-*$semicircuit in $\pi_\up$, which means 
that $\mu( E^+_{\up})=0$. As $\mu(E^+_{\tilde\pi})$ is either 0 or 1, 
the lemma follows. $\Box$

\section{Uniqueness of semi-infinite clusters}

Our next subject is the uniqueness of infinite 
clusters in half-planes, together with the stronger property that 
such clusters touch the boundary line infinitely often.
This was already a key result of Russo \cite{Rus}.
\begin{lem} {\bf (Line touching lemma)}\label{line_touching}
For any half-plane $\pi$, there exists $\G$-almost surely at most 
one infinite $+$ (resp.\ $+*$)
cluster $I^+_\pi$  (resp.\ $I^{+*}_\pi$) in $\pi$.
When it exists, this infinite cluster $\G$-almost surely
intersects the boundary line $\ell$ of $\pi$ infinitely often, in the 
sense that outside any finite $\D$ one can find an infinite path in 
this cluster starting from $\ell$.
\end{lem}
Just as Russo did, we derive this lemma from the absence of 
percolation for the $+$phase in the upper half-plane $\pi_\up$ with 
$-$boundary condition in $\pi_\up^c$ (which implies the uniqueness of 
the semi-infinite Gibbs measure, by the argument of Lemma \ref{lem:CNPR}). 
But for the latter we will give here a different argument using 
stochastic domination by a translation invariant Gibbs measure and 
Corollary \ref{cor:MM}. To state the result we write
$\pm$ for the configuration which is $+1$ 
on $\pi_\up$ and $-1$ on $\pi_\up^c$, and consider 
the semi-infinite limit
\be{pm}\mu^\pm_{\up}=\lim_{\D\uparrow\pi_\up}\mu_\D^{\pm}\ee 
which exists by stochastic monotonicity. 
\begin{lem}\label{semi-unique}
{\bf (No percolation on a bordered half-plane)}
$\mu^\pm_{\up}(E^{+*}_\up)=0$.
\end{lem}
{\sl Proof}. To begin we note that $\mu^\pm_{\up}$ is invariant 
under horizontal translations and stochastically maximal 
in the set of all Gibbs measures on $\pi_\up$ with $-$boundary 
condition in $\pi_\up^c$. This follows just as in the case of 
the plus-phase $\mu^+$ on the whole lattice. In particular,
$\mu^\pm_{\up}$ is trivial on the $\pi_\up$-tail $\T_\up=
\bigcap \{\F_{\pi_\up\setminus\L}:
{\L\subset\pi_\up\mbox{  finite}}\}$.
 We think of $\mu^\pm_{\up}$ as
a probability measure on $\O$ for which almost all configurations are 
identically equal to $-1$ on $\pi_\up^c$.

Next we consider the downwards translates 
$\mu^+_{n,-}=\mu^\pm_{\up}\circ \th_\v^{-n}$, $n\geq 0$.
Evidently, $\mu^+_{n,-}$ is obtained by an analogous infinite-volume 
limit in the half-plane $\{x_2\geq-n\}$. This shows that 
$\mu^+_{n,-}\preceq \mu^+_{n+1,-}$ by stochastic monotonicity,
so that the stochastically increasing limit 
$\mu^+_{-}=\lim_{n\to\infty}\mu^+_{n,-}$ exists.
Clearly $\mu^+_{-}\in\G$. Also, $\mu^+_{-}$ 
inherits the horizontal invariance of the $\mu^+_{n,-}$ and is in 
addition vertically invariant. Corollary \ref{cor:MM} therefore 
implies that $\mu^+_{-}=a\,\mu^- +(1-a)\mu^+$ for some coefficient 
$a\in[0,1]$.

We claim that $a>0$. For $n\ge1$ let $B_n$ denote the event that the 
origin is $-*$connected to the horizontal line $\{x_2=-n\}$. By the 
finite energy property and the
horizontal ergodicity of $\mu^+_{n,-}$, there exist for 
$\mu^+_{n,-}$-almost all $\o$ some random integers 
$m_\le(\o)<0<m_\ri(\o)$ 
such that $\o\equiv-1$ on 
$$
\s(\o)=\bigg\{x\in\Z^2: x_1\in\{m_\le(\o),m_\ri(\o)\},\;
-n\leq x_2\leq 0\bigg\}\;.
$$
Together with a segment of the line $\{x_2=-n-1\}$ on which $\o=-1$ 
$\mu^+_{n,-}$-almost surely, $\s(\o)$ forms a $-$semicircuit in 
$\pi_\lo$ surrounding the origin. An immediate application of the 
strong Markov property (applied to the largest such $\s$ in a large 
box) and the point-to-semicircuit lemma thus 
implies that $\mu^+_{n,-}(B_n)\geq 
\theta/2$. Therefore, writing $E^{-*}_{0,m}$ for the event that the 
origin belongs to some $-*$ cluster of size at least $m$ we find
$\mu^+_{n,-}(E^{-*}_{0,m})\geq \theta/2$ when $n\geq m$. Letting 
first 
$n\to\infty$ and then $m\to\infty$ we see that $\mu^+_{-}(E^{-*}) 
\geq \theta/2$.
Since $\mu^+(E^{-*})=0$ by Corollary \ref{ocean}, it follows that 
$a\geq\theta/2$, and the claim is proved.

To conclude the proof we observe that
\[
\mu^\pm_{\up}(E^{+*}_\up) \leq \mu^+_{-}(E^{+*}) = 1-a <1\;,
\]
again by Corollary \ref{ocean}.
Since $\mu^\pm_{\up}$ is trivial on $\T_\up$, the lemma follows.
$\Box$

\medskip\noindent
We are now able to prove Lemma \ref{line_touching} along the lines of 
Russo \cite{Rus}.

\medskip\noindent
{\sl Proof of Lemma \ref{line_touching}}. For definiteness we assume that 
$\pi=\pi_\up$; other half-planes merely correspond to a change of 
coordinates. We consider only infinite 
$+$clusters in $\pi_\up$; the case of $+*$clusters  
is similar. It is also clear that any result proved for the 
$+$sign is also valid with the $-$sign.

{\em Uniqueness: }The uniqueness of infinite $+$clusters in 
$\pi_\up$ is a consequence of the second statement, the line-touching 
property for infinite $-*$clusters. Indeed,
suppose there exists no infinite $-*$cluster in $\pi_\up$; then each 
finite set in $\pi_\up$ is surrounded by a $+$semicircuit, so that 
any two infinite $+$paths are necessarily $+$connected to each other. 
In the alternative case when an infinite $-*$cluster $I^{-*}_\up$ in 
$\pi_\up$ exists, 
this $I^{-*}_\up$ meets $\ell_\le$ or $\ell_\ri$ infinitely often, so 
that each infinite $+$cluster must meet the other half-line 
infinitely often. Hence, two such $+$clusters must cross each 
other, and are thus identical. 

{\em Line touching: } Let $\mu\in\Gex$ and $x\in\pi_\up$ and 
consider the event $A^+_x$ that $x$ 
belongs to an infinite $+$cluster in $\pi_\up$ which does not touch 
the horizontal axis $\ell_\h$. We will show that $\mu(A^+_x)=0$. 
Once this is 
established, we can take the union over all $x$ and use the finite 
energy property to see that for each finite $\D$ the event 
``an infinite $+$cluster in $\pi_\up$ is not connected to $\ell_\h$ 
outside $\D$'' also has probability zero, which means that almost 
surely any infinite 
$+$cluster in $\pi_\up$ must meet $\ell_\h$ infinitely often.

Intuitively, if $A^+_x$ occurs then the infinite $+$cluster 
containing $x$ is separated from $\ell_\h$ by an infinite $-*$path; 
but 
the spins `above' this path feel only the $-$boundary condition and 
thus believe to be in the $-$phase $\mu^-$, so that they will not 
form an infinite $+$cluster.

To make this intuition precise we fix some integer $k\geq1$ and 
consider the event $A^+_{x,k}$ that $x$ belongs to a 
$+$cluster of size at least $k$ which does not meet $\ell_\h$. Take a 
box $\D\subset\pi_\up$ containing $x$ and so large that there exists 
no path of length $k$ from $x$ to $\D^c$.
For $\o\in A^+_{x,k}$ we consider the largest set $\Ga(\o)\subset 
\D$ containing $x$ such that $\o=-1$ on 
$\partial\Ga(\o)\setminus\partial_\up\D$, where 
$\partial_\up\D=\partial\D\cap \pi_\up$. 
We also consider the event $E^+_{x,k}$ that $x$ belongs to a $+$cluster in 
$\pi_\up$ of size at least $k$. Using the fact that $A^+_{x,k}$ is
contained in the $\F_\D$-measurable event $\{\Ga\mbox{ exists}\}\cap
E^+_{x,k}$, we obtain by
the strong Markov property and the stochastic monotonicity of Gibbs 
distributions that
\[
\mu(A^+_{x,k})\leq \mu\bigg(\mu^\cdot_\Ga(E^+_{x,k}) \bigg)
\leq \mu_\D^{\pm}(E^+_{x,k})\;,
\]
where the $\pm$ boundary condition is defined as in \rf{pm}.
 Now, taking first the limit 
$\D\uparrow\pi_\up$ as in \rf{pm} and then letting $k\to\infty$ we find 
that $\mu(A^+_{x})\leq \mu^\pm_{\up}(E^+_\up)$. But the last 
expression vanishes by Lemma \ref{semi-unique}. $\Box$

\medskip\noindent
The butterfly lemma and shift lemma together still leave the 
possibility that all infinite butterflies have the same orientation, 
either horizontal or vertical. As a consequence of the line touching 
lemma, we can now show that both orientations must occur. 
\begin{lem}\label{orthogonal-butterflies} 
{\bf(Orthogonal butterflies)}
$\G$-almost surely there exist both a horizontal infinite butterfly 
in $\pi_\up$ and $\pi_\lo$ as well as a vertical infinite butterfly 
in $\pi_\le$ and $\pi_\ri$.
\end{lem}
{\sl Proof}. Suppose there exists some $\mu \in \Gex$ having almost 
surely no vertical infinite butterfly.
By the first step in the proof of the butterfly lemma, it then follows 
that $\mu=\mu\circ R_{k,\v}\circ T$ for all $k\in\Z$, and thus
$\mu=\mu\circ\th_{\h}^{-2}$.
By the tail triviality, $\mu$ is in fact ergodic under $\th_{\h}^2$; 
cf.\ Proposition (14.9) of \cite{Gii}.
By the butterfly lemma, horizontal infinite
butterflies do exist, say of color $+$.

We now use an argument similar to that in Corollary \ref{cor:MM},  
with the line touching lemma in place of the Burton-Keane theorem.
Fix any $n\ge1$. For $k\in\Z$ let $A_k$	denote the event that all 
spins along the straight path $p_{k,n}=((k,l):l = -n,\ldots, n)$ 
are $+1$, $(k,n)$ belongs to an infinite $+$cluster in
$\pi_{n,\up}=\{x\in\Z^2:x_2\ge n\}$,  and  $(k,-n)$ belongs to an 
infinite $+$cluster in $\pi_{n,\lo}=\{x_2\leq -n\}$.
Let $A$ be the event that $A_k$ occurs for infinitely many $k<0$ 
and infinitely many $k>0$. The finite energy property then shows 
that $\mu(A_0)>0$, and the horizontal ergodicity and Poincar\'e's
recurrence theorem (or the ergodic theorem) imply that
$\mu(A)=1$. But the line 
touching lemma guarantees that the infinitely many doubly-infinite 
`vertical' $+$paths passing through the horizontal axis are connected 
to each other in $\pi_{n,\up}$ and $\pi_{n,\lo}$. As	
$n$ was	arbitrary, it follows that 
almost surely each finite set is surrounded by a  
$+$circuit, and an infinite $-$cluster cannot exist. 
In view of Lemma \ref{lem:CNPR}, this implies that
$\mu=\mu^+$. But $\mu^+$ is not invariant under $R_\v\circ T$, 
in contradiction to what we derived for $\mu$. $\Box$

\medskip\noindent
The preceding argument can be used to derive the result of
Russo \cite{Rus} that $\mu^+$ and $\mu^-$ are the only phases which 
are periodic in one direction. We will not need this intermediate 
result.

\section{Non-coexistence of phases}
\label{inv}

In this section we will prove the following proposition.
\begin{prop} \label{invariance}
{\bf (Absence of non-periodic phases) }
Any Gibbs measure $\mu\in\G$ is invariant under translations, i.e., 
$\mu=\mu\circ\th_\h^{-1}$ and $\mu=\mu\circ\th_\v^{-1}$.
\end{prop}
Together with Corollary \ref{cor:MM} this will immediately imply the 
main theorem that each Gibbs measure is a mixture of the two phases 
$\mu^+$ and $\mu^-$. Our starting point is the following lemma 
estimating the probability that a semi-infinite cluster can be pinned 
at a prescribed point.
\begin{lem}
{\bf (Pinning lemma) } \label{pinning} 
Let $\mu\in\G$, and suppose that $\mu$-almost surely 
there exists an infinite $+*$cluster 
$I^{+*}_\up$ in $\pi_\up$ which meets the right semiaxis $\ell_\ri$ 
infinitely often. Then for each finite
square $\D=[-n,n]^2$ and $x\in\ell\ri$ we have
\[
\mu\bigg(\mbox{\rm $x$ is $+*$connected in 
$(\D\cup\ell_\le)^c$ to $I^{+*}_\up$}\bigg)
\ge \theta/4
\]
provided $x$ lies sufficiently far to the right. The same 
holds when `left' and `right' or `up' and `down' are interchanged.
\end{lem}
{\sl Proof}. By hypothesis, the infinite component of 
$I^{+*}_\up\setminus\D$ almost surely contains infinitely many points 
of $\ell_\ri$.
Thus, if $x\in\ell_\ri$ is located far enough to the right 
then, with probability exceeding $1/2$, at least one such 
point can be found left from $x$, and another such point can be found 
right from $x$. This means that $x$ is surrounded by a 
$+*$semicircuit $\s$ in $\pi_\up$
which belongs to $I^{+*}_\up$ and satisfies 
$\D\cap\mbox{Int\,}\s=\emptyset$.

Let $\L$ be a large square box containing $x$. If $\L$ is large 
enough, a semicircuit $\s$ as above can be found within $\L$ with 
probability still larger than $1/2$. We then can assume that $\s$ 
has the largest interior among all such semicircuits in $\L$. Using 
the strong Markov property and the point-to-semicircuit lemma we get 
the result.
$\Box$

\medskip\noindent
Our main task in the following is to analyze the situation when a 
half-plane contains both an infinite $+$cluster and an infinite 
$-$cluster.
(The line-touching lemma still allows this possibility.) 
In this situation it is useful to consider contours.

As is usually done in the Ising model, we draw lines of unit 
length between adjacent spins of 
opposite sign. We then obtain a system of polygonal curves running 
through the sites of 
the dual lattice $\Z^2+(\frac12,\frac12)$. A {\em contour\/} $\g$ in 
the upper half-plane $\pi_\up$ is a part of 
these polygonal curves which separates a $-$cluster in $\pi_\up$ from a 
$+*$cluster in $\pi_\up$. This corresponds to the convention that at 
crossing points the contours 
are supposed to bend around the $-$spins. (The artificial asymmetry 
between $+$ and $-$ does not matter, and we could clearly make the 
opposite convention.) On its two sides, $\g$ is accompanied by a 
$+*$quasipath $f^{+}_\g$ and a $-$quasipath $f^{-}_\g$ which will be 
called the $+$ resp.\ $-$face of $\g$; the prefix `quasi' 
indicates that the faces are not necessarily self-avoiding but may 
contain loops.
\begin{lem}\label{contours}
{\bf (Semi-infinite contours) } $\G$-almost surely on $E^{+*}_\up 
\cap E^-_\up$ there exists a unique semi-infinite contour $\g_\up$ 
in $\pi_\up$. $\g_\up$ starts between two points of the 
horizontal axis $\ell_\h$ and intersects each horizontal line in 
$\pi_\up$ only finitely often.
\end{lem}
{\sl Proof}. Let $I^{+*}_\up$ be the unique infinite $+*$cluster in 
$\pi_\up$, and $I^-_\up$ the unique infinite $-$cluster in $\pi_\up$.
For definiteness we assume that $I^{+*}_\up$ meets $\ell_\ri$ 
infinitely often, and $I^-_\up$ meets $\ell_\le$ infinitely often.
Let $x$ be the rightmost point of $I^-_\up\cap\ell_\h$ and $\g_\up$ the 
contour in $\pi_\up$ starting from the line segment separating $x$ and 
$y=x+(1,0)$. Since $I^-_\up$ contains an infinite $-$path starting 
from $x$ which cannot be traversed by $\g_\up$, $\g_\up$ cannot return to 
$\ell_\h$ on the left-hand side of $x$. But $\g_\up$ can also not return 
to $\ell_\h$ on the right-hand side of $y$, since otherwise the 
$-$face of $\g_\up$ would establish a $-$connection in $I^-_\up$ from $x$ to 
a point of $\ell_\h$ to the right of $y$, in contradiction to the 
choice of $x$. Hence $\g_\up$ can never end and must therefore be infinite.

Let $\g$ be any infinite contour in $\pi_\up$. Then 
the infinite $-$face $f^-_{\g}$ must belong to $I^-_\up$,
by the uniqueness of the infinite $-$cluster. This implies that $f^-_{\g}$ 
must lie on the ``left-hand 
side'' of $\g_\up$. Likewise, the $+*$face $f^{+*}_{\g}$ must belong 
to the ``side on the right'' of $\g_\up$. Hence $\g=\g_\up$, proving the 
uniqueness of $\g_\up$. 

Finally, let $\ell=\{x_2=n\}$, $n\geq 1$, be a horizontal line in 
$\pi_\up$ and $\pi=\{x_2\geq n\}$ the half-plane above $\ell$. By the 
shift lemma and the above, $\pi$ contains a unique semi-infinite 
contour $\g$ starting from the line segment between two adjacent 
points $u$ and $v$ of $\ell$. $u$ and $v$ belong to the infinite faces 
of $\g$ and therefore to $I^{+*}_\up$ resp.\ $I^-_\up$. By the line 
touching lemma, this means that $u$ and $v$ are $+*$connected resp. 
$-$connected to the axis $\ell_\h$. The unique continuation of $\g$ can 
therefore visit only finitely many sites of $\pi_\up$, 
and thus must reach $\ell_\h$ after finitely many steps; this 
continuation is then equal to $\g_\up$, by the uniqueness of the 
latter. This shows that $\g_\up$ visits the line $\ell$ only finitely 
often. $\Box$

\medskip\noindent
From now on we consider a fixed extremal Gibbs measure $\mu\in\Gex$.
We want to prove that $\mu$ is horizontally invariant. (The proof of 
vertical invariance is similar.) To this end we consider its 
horizontal translate
$\hat\mu=\mu\circ\th_\h^{-1}$, as well as the product measure
$\hat\nu=\mu\otimes\hat\mu$ on $\O\times\O$. It is convenient to 
think of the latter as a duplicated system consisting of two 
independent layers.
The following lemma is a slight variation of a result of Aizenman 
\cite{Aiz} in his proof of the main theorem; our proof differs in part.
\begin{lem}\label{intersections}
{\bf (Fluctuations of the semi-infinite contour) }
Suppose $\pi_\up$ contains a semi-infinite contour $\g_\up$ 
$\mu$-almost surely. Then
for $\hat\nu$-almost all 
$(\o,\hat\o)\in\O^2$, $\g_\up(\o)$ and $\g_\up(\hat\o)$ intersect 
each other infinitely often.
\end{lem}
{\sl Proof}. By tail triviality, we can assume that $\g_\up$ 
has its $+$face on the left-hand side almost surely; the alternative 
case is analogous. For any $n\geq 1$ we let 
\[
a_n = \max\{k\in\Z: (k,n)\in I^{+*}_{\pi_{n,\up}}\}
\]
be the abscissa of the point at which $\g_\up$ enters definitely into 
the half-plane $\pi_{n,\up}=\th_\v^n\pi_\up$ above the height 
$n$. Consider the product measure $\nu=\mu\otimes\mu$ and the event
\[
F =\{(\o,\o')\in\O^2: \mbox{$\g_\up(\o)$ and $\g_\up(\th_\h\,\o')$ meet 
each other only finitely often} \}\;.
\]
We need to show that $\nu(F)=0$. 

Suppose that $F$ occurs. Then $\g_\up(\o)$ lies strictly on one side 
of $\g_\up(\th_\h\,\o')$ above some random level $n$. 
Hence we have either $a_n(\o) > a_n(\th_\h\,\o')$
eventually, or $a_n(\o) < a_n(\th_\h\,\o')$ eventually. Using the 
abbreviation $d_n(\o,\o')=a_n(\o)-a_n(\o') = 
a_n(\o)-a_n(\th_\h\,\o')+1$, we thus see that
\[
F \subset A\cup B \equiv \{ d_n \geq 0 \mbox{ eventually}\} \cup 
\{ d_n \leq 0 \mbox{ eventually}\}\;.
\]

Suppose now that $\nu(F)>0$. Then, by symmetry, $\nu(A)=\nu(B)>0$. By 
the tail-triviality of $\mu$, it follows that $\nu(A)=\nu(B)=1$. 
This is because $A,B$ are measurable with respect 
to the `product-tail' $\T^{(2)}=\bigcap\{\F_{\L^c}\otimes\F_{\L^c}
:\L\subset\Z^2\mbox{ finite}\}$ in $\O^2$, which is trivial by 
Fubini's theorem. (One should not be mistaken to 
believe that $A$ was measurable with respect to the smaller 
`tail-product' $\T\otimes\T$. It is only the case that the 
$\o$-section $A_\o$ of $A$ belongs to $\T$ for any $\o$, and the 
function $\o\to\mu(A_\o)$ is $\T$-measurable.) We thus conclude 
that $\nu(A\cap B)=1$, meaning that $d_n=0$ eventually almost surely.
The lemma will therefore be proved once we have shown that this is 
impossible.

To this end we claim first that $\nu(d_{n}\geq1)\geq \d\,\nu(d_{n+1}=0)$ 
for all $n$ and some constant $\d>0$.
To see this let $A_{k,n}=\{(\o,\o'):a_{n+1}(\o)=a_{n+1}(\o')=k \}$,
$\D_{k,n}$ the two-point set consisting of the points $(k,n)$ and 
$(k+1,n)$, and $B_{k,n}$ the event that $\o=(+1,+1)$ on $\D_{k,n}$ 
and $\o'=(+1,-1)$ on $\D_{k,n}$; see the figure. 
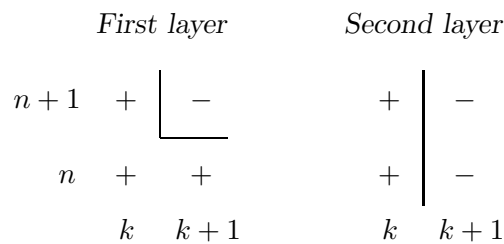
\setlength{\unitlength}{5mm}
\[
\begin{picture}(11,6.5)
\put(1,5.5){{\sl First layer}}
\put(7.6,5.5){{\sl Second layer}} 
\put(0,1.5){$n$} \put(-1.2,3.5){$n+1$}
\put(1.6,0){$k$} \put(3.1,0){$k+1$}
\put(8.6,0){$k$} \put(10.1,0){$k+1$}
\multiput(1.5,1.5)(0,2){2}{$+$}
\put(3.5,1.5){$+$} \put(3.5,3.5){$-$}
\multiput(8.5,1.5)(0,2){2}{$+$}
\multiput(10.5,1.5)(0,2){2}{$-$}
\put(2.7,2.7){\line(0,1){1.8}}
\put(2.7,2.7){\line(1,0){1.8}}
\put(9.7,0.9){\line(0,1){3.6}}
\end{picture}
\]
We then have
\[
\nu(B_{k,n}|\F_{\D_{k,n}^c}\otimes\F_{\D_{k,n}^c})(\o,\o')=
\mu^\o_{\D_{k,n}}\otimes \mu^{\o'}_{\D_{k,n}}(B_{k,n})\geq 
[1+e^{8\b}]^{-4}\equiv\d
\]
and thus
\[
\nu(\{d_{n}\geq1\}\cap A_{k,n})
\geq \nu\Big( \nu(B_{k,n}|\F_{\D_{k,n}^c}\otimes\F_{\D_{k,n}^c})\; 1_{A_{k,n}}
\Big)\geq\d \; \nu(A_{k,n})
\]
because $A_{k,n}$ is an event in $\D_{k,n}^c$. Summing over $k$ we get the 
claim. 

Now, if $d_n=0$ eventually almost surely then
\[
\liminf_{n\to\infty} \nu( d_n\ge 1) \ge \d\, \liminf_{n\to\infty}\nu( d_{n+1}=0) = 
\d\;,
\]
so that with positive probability we have simultaneously $d_n\geq 1$ 
infinitely often and $d_n=0$ eventually. Since this is impossible,
we conclude that $\nu(F)=0$.
$\Box$

\medskip\noindent 
The following percolation result for the 
duplicated system with distribution $\hat\nu$ was already a 
cornerstone of Aizenman's argument \cite{Aiz}. We prove it here 
differently, avoiding his use of the fact that the limiting Gibbs measure 
for the $\pm$boundary condition is translation invariant.
We will say that a path in $\Z^2$ is a $\ls$path 
for a pair $(\o,\hat\o)\in\O^2$ if $\o(x)\leq\hat\o(x)$ for all its 
sites $x$. In the same way we define $\lss$paths, and we can speak 
of $\lss$circuits and $\lss$clusters. 
\begin{lem}\label{lss-circuits}
{\bf(No $(+,-)$percolation in the duplicated system) } 
$\hat\nu$-almost surely each finite square $\D=[-n,n]^2$ is 
surrounded by 
a $\lss$circuit in $\Z^2$. 
\end{lem}
{\sl Proof}. Consider any two points $x\in\ell_\le$ and 
$y\in\ell_\ri$.
We claim that with $\hat\nu$-probability at least $(\theta/4)^2$ 
there exists a $\lss$path from $x$ to $y$ `above' $\D$, provided $x$ 
and $y$ are located sufficiently far to the left resp.\ to the right. 
We distinguish three cases. 

{\em Case 1: $\mu(E^+_\up)=0$}. By Lemma 
\ref{orthogonal-butterflies}, $\pi_\up$ then almost surely contains an 
infinite $-$cluster $I^-_\up$, and each finite subset of $\pi_\up$ is 
surrounded by a $-*$semicircuit in $\pi_\up$. In other words, an 
infinite $-*$cluster $I^{-*}_\up$ in $\pi_\up$ exists and touches 
both $\ell_\le$ and $\ell_\ri$ infinitely often. By the pinning lemma 
and the positive correlations of $\mu$,
with $\mu$-probability at least $(\theta/4)^2$ both $x$ and $y$ are 
$-*$connected to $I^{-*}_\up$ outside $\D$, and therefore also 
$-*$connected to each other by a $-*$path $p$ above $\D$. However, 
this $-*$path $p$ on the first layer is certainly also a $\lss$path 
for the duplicated system, and the claim follows.

{\em Case 2: $\mu(E^-_\up)=0$}. In this case we also have
$\hat\mu(E^-_\up)=0$. Interchanging $+$ and $-$ and replacing $\mu$ 
by $\hat\mu$ in Case 1, we find that with $\hat\mu$-probability at 
least $(\theta/4)^2$, there exists a $+*$path $\hat p$ in the second 
layer above $\D$ from $x$ to $y$. Since $\hat p$ 
is again a $\lss$path 
for the duplicated system, the claim follows as in the first case.

{\em Case 3: $\mu(E^+_\up)=\mu(E^-_\up)=1$ }. 
Then $\mu$-almost surely there exists a unique semi-infinite contour 
$\g_\up$, and by tail triviality we can assume (for definiteness) 
that $\g_\up$ has its $+$face on the left-hand side $\mu$-almost 
surely, and thus also $\hat\mu$-almost surely. By the pinning lemma 
and the independence of the two layers, the following event has 
$\hat\nu$-probability at least $(\theta/4)^2$: 
\bit
\item[--]in the first layer, $y$ is $-*$connected off $\D$ to 
$I^{-}_\up(\o)$, and thus to the $-$face $f^-_\up(\o)$ of 
$\g_\up(\o)$; that is, there exists an infinite $-*$path $p_y^-(\o)$ 
from $y$ outside $\D$ eventually running along $\g_\up(\o)$;
\item[--]in the second layer, $x$ is $+*$connected off $\D$ to 
$I^{+*}_\up(\hat\o)$, and thus to the $+$face $f^+_\up(\hat\o)$ of 
$\g_\up(\hat\o)$; that is, there exists an infinite $+*$path 
$p_x^+(\hat\o)$ from $x$ outside $\D$ eventually running along 
$\g_\up(\hat\o)$.
\eit
Since $\g_\up(\o)$ and $\g_\up(\hat\o)$ intersect each other 
infinitely often by Lemma \ref{intersections}, the union of 
$p_y^-(\o)$ and $p_x^+(\hat\o)$ contains a 
$*$path from $x$ to $y$ which by construction is a $\lss$path for the 
duplicated system. This proves the claim in the final case. 

To conclude the proof of the lemma, we let $A_{x,y}$ denote the event 
that there exist a $\lss$path from $x$ to $y$ above $\D$, and 
$B_{x,y}$ the event that such a path exists below $\D$. The indicator 
functions of these events can be written as increasing functions $f$ 
resp.\ $g$ of the difference configuration $\hat\o-\o$. Using the 
positive correlations of $\hat\mu$ and $\mu$ we thus obtain
\bea
\hat\nu(A_{x,y}\cap B_{x,y}) &=& \int\mu(d\o)\int\hat\mu(d\hat\o)\; 
f(\hat\o-\o)\, g(\hat\o-\o)\\
&\geq& \int\mu(d\o)\; \hat\mu(f(\cdot-\o))\,\hat\mu(g(\cdot-\o))\\
&\geq& \hat\nu(A_{x,y})\,\hat\nu(B_{x,y})
\ \geq\ (\theta/4)^4\;.
\eea

The last inequality follows from the claim and its analogue for the 
lower half-plane. However, if $A_{x,y}\cap B_{x,y}$ occurs then $\D$ 
is surrounded by a $\lss$circuit for the duplicated system. Letting 
$\D\uparrow\Z^2$ we see that with probability at least $(\theta/4)^4$ 
each finite set is surrounded by a $\lss$circuit. Since this event is 
measurable with respect to the product-tail $\T^{(2)}$ on which 
$\hat\nu$ is trivial, the lemma follows.
$\Box$

\medskip\noindent
It is now easy to complete the proof of Proposition \ref{invariance}
as in \cite{Aiz}.

\medskip\noindent
{\sl Proof of Proposition \ref{invariance}}. Consider any square 
$\D=[-n,n]^2$, and let $\e>0$. By Lemma \ref{lss-circuits}, $\D$ is 
$\hat\nu$-almost surely surrounded by a $\lss$circuit, and with 
probability at least $1-\e$ such a $\lss$circuit can be 
found in a sufficiently large square $\L$. Let $\Ga$ be the interior 
of the largest such $\lss$circuit; if no such $\lss$circuit exists 
let $\Ga=\emptyset$. Then we find for any increasing 
$\F_\D$-measurable function $0\leq f\leq1$, using the strong 
Markov property of $\hat\nu$ and 
the fact that $\mu_\Ga^\o\preceq \mu_\Ga^{\hat\o}$ when 
$\Ga(\o,\o')\ne\emptyset$, 
\bea
\mu(f)&=&\hat\nu(f\otimes 1)\ \leq\ 
\int_{\{\Ga\ne\emptyset\}} d\hat\nu(\o,\o')\;
\mu^{\o}_{\Ga(\o,\o')}(f) +\e\\
&\leq&\int d\hat\nu(\o,\o')\;\mu^{\o'}_{\Ga(\o,\o')}(f) +\e
\ =\ \hat\nu(1\otimes f)+\e\ = \ \hat\mu(f)+\e\;.
\eea
Letting $\e\to0$ and $\D\uparrow\Z^2$ we find that 
$\mu\preceq\hat\mu$. Interchanging $\mu$ and $\hat\mu$ 
(i.e., the roles of the layers) we get the reverse relation. 
Hence $\mu=\hat\mu$, so that $\mu$ is horizontally invariant.
The vertical invariance follows similarly by an interchange of 
coordinates. $\Box$

\section{Extensions}
\label{extensions}

Which properties of the square lattice $\Z^2$ entered into 
the preceding arguments? The only essential feature was its 
invariance 
under the reflections in all horizontal and vertical lines with 
integer coordinates. 
We claim that the theorem remains true for the Ising model 
on any connected graph $\LL$ with these properties. (The Ising model 
on the triangular and hexagonal lattices has already been treated in 
\cite{Fuk}.)

To be more precise, let $\R=\{R_{k,\h},R_{k,\v}:k\in\Z\}$ denote the 
set of all reflections of the Euclidean 
plane ${\bf R}^2$ in horizontal or vertical lines with integer 
coordinates, and suppose $\LL$ is a locally finite subset 
of ${\bf R}^2$ 
which (after suitable scaling and rotation) is $R$-invariant
for all $R\in\R$.
Such an $\LL$ is uniquely determined by its finite intersection 
 with the unit cube $[0,1]^2$, and it is periodic with period $2$. 
Suppose further that $\LL$ is equipped with a symmetric 
neighbor relation `$\sim$' satisfying
\bit
\item[(L1)] each $x\in\LL$ has only finitely many `neighbors' 
$y\in\LL$ 
satisfying $x\sim y$;
\item[(L2)] $x\sim y$ if and only if $Rx\sim Ry$ for all $R\in\R$;
\item[(L3)] $(\LL,\sim)$ is a connected graph.
\eit
If $x\sim y$ we say that $x$ and $y$ are connected by an edge, which 
is visualized by the straight line segment between $x$ and $y$.
The preceding assumptions simply mean that $(\LL,\sim)$ is a locally 
finite connected graph admitting the reflections $R\in\R$, and thereby 
the translations $\th_x$, $x\in 2\Z^2$, as graph automorphisms. 
The fundamental further assumption is 
\bit
\item[(L4)] $(\LL,\sim)$ is planar, i.e., the edges in ${\bf R}^2$ between 
different pairs of neighboring points have only endpoints in common.
\eit
The complement (in ${\bf R}^2$) of 
the union of all edges then splits into connected components called 
the faces of $(\LL,\sim)$.

As will be explained in more detail in the appendix, the properties 
(L1) to (L4)
are sufficient for all geometric arguments above. Some particular
examples are
\medskip

$\bullet$ the {\em triangular lattice\/} {\bf T}. This is the 
$\R$-invariant lattice satisfying 
${\bf T}\cap[0,1]^2=\{(1,0),(0,1)\}$  and 
$(-1,0)\sim(1,0)\sim(0,1)\sim(2,1)$; the 
remaining edges result from (L2). 

$\bullet$ the hexagonal or {\em honeycomb lattice\/} {\bf H}. Here,
for example, ${\bf H}\cap[0,1]^2=\{(\frac13,1),(\frac23,0) \}$
and $(-\frac13,1)\sim(\frac13,1)\sim(\frac23,0)\sim
(\frac43,0)$; all other edges are 
again determined by (L2). 

$\bullet$ the {\em diced lattice}. This is obtained from the 
honeycomb lattice by placing points in the centers of the hexagonal 
faces and connecting them to the three points in the west, north-east 
and south-east of these faces; to obtain reflection symmetry an 
additional shift by $(-\frac13,0)$ is necessary.
See p.\ 16 of \cite{Kes} for more details.

$\bullet$ the covering lattice of the honeycomb lattice, the 
{\em Kagom\'e lattice}, cf.\ p.\ 37 of \cite{Kes}.

\bigskip\noindent
As for the interaction, 
it is neither necessary that all adjacent spins interact in the 
same way, nor that the interaction is invariant under the spin flip.
Except for attractivity,
we need only the invariance under {\em simultaneous\/} flip-reflections 
(which in particular implies periodicity with period 2).
As a result, 
we can consider any system of spins $\o(x)=\pm1$ with formal 
Hamiltonian of the form
\be{Ham}
H(\o)=\sum_{x\sim y} U_{x,y}(\o(x),\o(y)) +\sum_{x\in\LL} 
V_x(\o(x))\;,
\ee
where for all $a,b\in\{-1,1\}$ we have $U_{x,y}(a,b)=U_{y,x}(b,a)$ and
\bit
\item[(H1)] $U_{x,y}(1,\cdot)-U_{x,y}(-1,\cdot)$ is decreasing on 
$\{-1,1\}$;
\item[(H2)] $U_{x,y}(a,b)=U_{Rx,Ry}(-a,-b)$ and $V_x(a)=V_{Rx}(-a)$
for all $R\in\R$.
\eit
Assumption (H1) implies that the FKG inequality is applicable,
and (H2) expresses the invariance under simultaneous spatial 
reflection and spin flip.
We thus obtain the following general result.
\begin{thm}\label{general}
Consider a planar graph $(\LL,\sim)$ as above and an 
interaction of the form \rf{Ham} satisfying (H1) and (H2). 
Then there exist no more than two extremal Gibbs measures.
\end{thm}
The standard case, of course, is the ferromagnetic Ising model 
without external field; this corresponds to the choice 
$U_{x,y}(a,b)=-\b ab$ and $V_x\equiv 0$. 

But there is also 
another case of particular interest. Consider 
$\LL=\Z^2+(\frac12,\frac12)$, the shifted square lattice with its 
usual 
graph structure. $\LL$ is bipartite, in the sense that $\LL$ splits 
into two disjoint sublattices, $\LL_{even}$ and $\LL_{odd}$, such that 
all edges run from one sublattice to the other. If we set 
$U_{x,y}(a,b)=-\b ab$ 
and define a staggered external field
\[
V_x(a)=\left\{\ba{rl}-h a&\mbox{if }x\in \LL_{even}\\
h a& \mbox{if }x\in \LL_{odd}\ea \right.
\]
with $h\in{\bf R}$ then the conditions (H1) and (H2) hold; 
here we take advantage of the fact that the 
reflections $R\in\R$ map 
$\LL_{even}$ into $\LL_{odd}$ and vice versa. 
But it is well-known that this model is isomorphic to the 
{\em anti\/}ferromagnetic Ising model on $\Z^2$ with homogeneous 
external field $h$; the 
isomorphism consists in flipping all spins in $\LL_{odd}$. This 
gives us the following result.
\begin{cor}\label{AF}
For the Ising antiferromagnet on $\Z^2$ for any inverse temperature 
and arbitrary external field there exist at most two extremal 
Gibbs measures.
\end{cor}
This corollary does not extend to non-bipartite lattices such as the 
triangular lattice. In fact, for the Ising antiferromagnet on {\bf T} 
one expects the existence of three different phases for suitable $h$.

Another repulsive model to which our arguments can be applied is
the hard-core lattice gas on $\Z^2$, which is also known as the hard 
square model. In this model, the values 
$-1$ and $1$ are interpreted as the absence resp.\ presence of a 
particle, and no particles are allowed to sit on adjacent sites.
Its Hamiltonian is of the form \rf{Ham} with
\[
U_{x,y}(a,b)=\left\{\ba{cl}\infty &\mbox{if }a=b=1,\\
                              0   &\mbox{otherwise,}\ea\right.
\qquad
V_x(a)=\left\{\ba{cl}-\log\lambda &\mbox{if }a=1,\\
                              0   &\mbox{otherwise;}\ea\right.
\]
The parameter $\lambda>0$ is called the activity.
Interchanging the values $\pm 1$ on $\LL_{odd}$ we obtain an 
isomorphic attractive model to which our techniques can be applied,
although the interaction takes the value $+\infty$ so that the 
finite energy condition does not hold as it stands. 
However, there are still enough admissible configurations to 
satisfy all needs of the Burton-Keane theorem and our other 
applications of the finite energy property; more details
will be provided in the appendix. 
We therefore can state the following theorem.
\begin{thm}\label{hard-core}
For the hard-core lattice gas on $\Z^2$ at any activity $\lambda>0$ 
there exist at most two extremal Gibbs measures.
\end{thm}
\section{Appendix}

Here we explain in more detail how our
arguments can be extended to obtain Theorems 
\ref{general} and \ref{hard-core}.

\medskip\noindent
{\sl Comments on the proof of Theorem \ref{general}. }
(1) {\em $*$Connections and contours. }A 
basic consequence of the planarity assumption (L4) is that $(\LL,\sim)$ 
admits a conjugate matching graph $(\LL,\stackrel{*}{\sim})$.  As 
indicated by the 
notation, this conjugate graph has the same set of vertices, but the 
relation $x\stackrel{*}{\sim} y$ holds if either $x\sim y$ or $x$ and 
$y$ are distinct points (on the border) of the same face of $(\LL,\sim)$. 
(Note that this matching dual is in general not planar. 
An interesting exception is the triangular lattice {\bf T}, which
is self-matching.) 
The edges of $(\LL,\stackrel{*}{\sim})$ are then used to define the 
concept of $*$connectedness. The construction implies that
every path in $(\LL,\sim)$ is also a $*$path (i.e., a path in 
$(\LL,\stackrel{*}{\sim})$), and that the outer boundary of any 
cluster is a $*$path, and vice versa.
(The latter property holds for arbitrary matching pairs of graphs
as defined in Kesten \cite{Kes}, for example. However, we also
used repeatedly the former property which does not extend to general 
matching pairs. In particular, this means that our results do not 
apply to the Ising model on the matching conjugate of $\Z^2$ having
nearest-neighbor interactions {\em and\/} diagonal interactions.) 

Another consequence of planarity is that we can draw contours
separating clusters from $*$clusters. Such contours can either be 
visualized by broken lines passing through the edges of $(\LL,\sim)$,
or simply as a pair consisting of a quasipath and an adjacent 
$*$quasipath, namely the two faces of the contour.

(2) {\em Half-planes and boundary lines. }A 
half-plane $\pi$ in $\LL$ is still defined as
the intersection of $\LL$ with a set of the form $\{x\in {\bf 
R}^2:x_i\geq k\}$, $k\in\Z$, $i\in\{1,2\}$, or with $\leq$ 
instead of $\geq$. However, the `boundary line' $\ell$ is now in 
general not a straight line but rather the set $\ell=\{x\in\pi:x\sim 
y \mbox{ for some } y\notin\pi\}=\partial (\pi^c)$. In particular, $\ell$ is  
not necessarily a line of fixed points for the reflection $R\in\R$ 
mapping $\pi$ onto its conjugate halfplane $\pi'$. Rather, 
for each $x\in\ell$ we have either $Rx=x$ or $Rx\sim x$. 
For example, for $\LL={\bf T}$, the triangular lattice,
$\pi_\up$ and $\pi_\lo$ have a common straight boundary line, but 
the boundaries of 
$\pi_\ri$ and $\pi_\le$ are not straight; besides a common part on 
the vertical axis they also contain the adjacent points $(1,k)$ 
resp.\ $(-1,k)$, $k\in 2\Z$.
For the honeycomb lattice {\bf H}, $\pi_\up$ and
$\pi_\lo$ have again a common straight boundary line, but $\pi_\ri$ and 
$\pi_\le$ have no common points.

Nevertheless, it is easy to see that Lemma \ref{point-to-semicircuit} 
(and thus also Lemma \ref{pinning}) are still valid, and these are the 
only results in which fixed points of reflections show up. 
In all other places one has only to observe that the axes $\ell_\h$ 
and $\ell_\v$ get a different meaning according to which 
half-space is considered; so one has to distinguish 
between an `upper' horizontal axis $\ell_{\h,\up}$ (being the 
boundary `line' of $\pi_\up$) and a `lower' 
horizontal axis $\ell_{\h,\lo}$, and similarly 
between $\ell_{\v,\le}$ and $\ell_{\v,\ri}$. 

(3) {\em Construction of connections and paths. }At 
various places we needed to establish prescribed connections
or to construct specific paths.
For example, the key idea of Lemma 
\ref{shift-lemma} was to extend $-*$semicircuits in $\tilde\pi$ to the 
boundary line of $\pi$. In the present setup, this will in general 
require a finite $-$path rather than a single $-$spin, so that one has 
to adapt the definition of $A_n$ accordingly. 
In view of (L3) this is obviously possible, and one will only end up 
with a higher power of $\d$. Similarly, the $-$semicircuit $\s$ 
in the proof of Lemma \ref{semi-unique} has in general to be redefined 
using the geometry of $\LL$, and the same is the case for 
the points $(k,0)$ in the definition of $A_k$ in the proof of Corollary 
\ref{cor:MM}, the paths $p_{k,n}$ in the proof of Lemma 
\ref{orthogonal-butterflies}, and the sets $\D_{k,n}$  in the proof of 
Lemma \ref{intersections}; see also comment (5) below. 

(4) {\em Flip-reflection invariance. }In 
the standard Ising model on $\Z^2$ it is true that the phases 
$\mu^+$ and $\mu^-$ are invariant under all $R\in\R$ and related to 
each other by the spin flip $T$. However, we did not make use of this fact, 
cf.\ the comments after Corollary 
\ref{ocean}. We only needed that $\mu^+=\mu^-\circ 
R\circ T$ for all $R\in\R$ (implying that $\mu^+$ and $\mu^-$ are 
periodic, and that any flip-reflection invariant $\mu$ is different 
from these phases; the latter was used in Lemmas \ref{butterfly} and 
\ref{orthogonal-butterflies}). 
This, however, already holds whenever the interaction is only invariant
under simultaneous flip-reflections, as stated in assumption (H2). 
This property is also sufficient for 
flip-reflection domination and the point-to-semicircuit lemma, as their
proofs only use the composed mappings $R\circ T$ for $R\in\R$.

(5) {\em Translations. }Since the lattice and the 
interaction are in general only preserved by the translation subgroup
$\th_x$, $x\in 2\Z^2$, we have to confine ourselves to this class of 
translations. We did this already in the proof of the butterfly lemma and 
its Corollary \ref{cor:MM}, and we can obviously do so in the proof of
Lemma \ref{semi-unique}. The only statements needing discussion are
Proposition \ref{invariance} and Lemma \ref{intersections}. 
The former now only asserts that each Gibbs measure is periodic with 
period 2. Accordingly, in Lemma \ref{intersections} and below we have 
to replace $\th_\h$ by $\th_\h^2$. In addition, the minimal distance 
between distinct lattice points can be less than 1, and the origin 
does not necessarily belong to $\LL$. So, $a_n$ has to be defined as 
the abscissa of the rightmost point in the boundary line of 
$\pi_{n,\up}$ which belongs to $I_{\pi_{n,\up}}^{+*}$, and
$d_n(\o,\o')=a_n(\o)-a_n(\o') = a_n(\o)-a_n(\th_\h^2\,\o')+2$.
In general, we then have only the inclusion
\[
F \subset \{ d_n > -2 \mbox{ eventually}\} \cup 
\{ d_n < 2 \mbox{ eventually}\}\;,
\]
and we need to derive a contradiction from the assumption that 
$|d_n|< 2$ eventually almost surely. This means that we have to 
prescribe the configurations for the two layers on larger sets 
than $\D_{k,n}$ (depending on $n$ and both
$a_n(\o)$  and $a_n(\o')$) to obtain the inequality 
$\nu(|d_{n}|\geq 2)\geq \d\,\nu(|d_{n+1}|<2)$ for some $\d>0$.
While this is tedious to write down in full generality, it 
should be clear how it can be done.

\medskip\noindent
{\sl Comments on the proof of Theorem \ref{hard-core}. }Just
as in the case of the Ising antiferromagnet, we replace $\Z^2$
by its translate $\LL=\Z^2+(\frac12,\frac12)$. So we make sure
that all reflections $R\in\R$ map $\LL_{even}$ into $\LL_{odd}$ 
and vice versa. Nevertheless, below it will be convenient 
to ignore the shift by $(\frac12,\frac12)$ and to characterize the 
lattice points by integer coordinates.
Performing a spin flip on $\LL_{odd}$ we obtain an 
isomorphic model which is defined by setting
\[
U_{x,y}(a,b)=\left\{\ba{cl}\infty &\mbox{if }a=\epsilon(x), 
                                             b=\epsilon(y),\\
                              0   &\mbox{otherwise,}\ea\right.
\qquad
V_x(a)=\left\{\ba{cl}-\log\lambda &\mbox{if }a=\epsilon(x),\\
                              0   &\mbox{otherwise,}\ea\right.
\]
where $\epsilon(x)=1$ if $x\in\LL_{even}$ and $\epsilon(x)=-1$ 
otherwise.
This model satisfies both (H1) and (H2). However, the
finite energy condition is violated because  $U_{x,y}$
takes the value $+\infty$. Let us see how this obstacle
can be overcome. The basic observation is that the `vacuum 
configuration' $-\epsilon$ can occur in any finite region with 
positive probability.

(1) In the proof of 
the Burton-Keane theorem, the finite energy property is used to 
connect different $+$clusters with positive probability. This is 
still possible because for any box $\D$, 
any $x\in\D$, any finite 
number of points $x_1,\ldots,x_k\in\partial\D$, and any $\o$ with 
$\o(x_1)=\ldots =\o(x_k)=+1$ we have 
$$
\mu_\D^\o(x \mbox{ is $+$connected to }x_1,\ldots,x_k)>0\;.
$$

(2) A different use of the finite energy property is made in the proofs of
Corollary \ref{cor:MM} and Lemma \ref{orthogonal-butterflies}: 
the events $A_k$ there involve the existence of both $+$ and $-$paths.
To adapt the proof of Corollary \ref{cor:MM} to the present case we 
redefine $A_0$ as the event
that a prescribed point $x\in\LL_{odd}$ belongs to a two-sided 
infinite $+$path with its two halves staying in $\pi_\up$ resp.\ 
$\pi_\lo$, and a neighbor point $y\in\LL_{even}$ belongs 
to an infinite $-$cluster; for $k\in 2\Z$ we set $A_k=\th_\h^{-k} A_0$.
A $+$spin at $x$ then
does not interfere with a $-$spin at $y$. Therefore, if $\D$ is a
sufficiently large box and $u_1,u_2, u_3\in\partial\D$ are three points 
belonging to infinite $+$, $+$ resp.\ $-$clusters meeting $\D$, 
we can find paths $p_1,p_2$ in $\D$ from $x$ to 
$u_1$ resp.\ $u_2$ and a path $p_3$ from $y$ to $u_3$ such that 
$y$ is the only site of $p_3$ which is adjacent to $p_1\cup p_2$.
The 
configuration in $\D$ which is equal to $+1$ on $p_1\cup p_2$, $-1$ 
on $p_3$, and $-\epsilon$ otherwise then has positive conditional 
probability given the configuration in $\D^c$. 
This shows that $\mu(A_0)>0$. The proof of Lemma \ref{orthogonal-butterflies} 
can be adapted in a similar manner.

(3) In Lemma \ref{semi-unique} we used the finite energy property to 
make sure that $\mu_{n,-}^+(\o\equiv -1\mbox{ on }p)>0$, where 
$p=\{0\}\times\{-n+1,\ldots,0\}$. To obtain the same result here we 
simply set $\D=p\cup\partial p\setminus\{x_2=-n\}$ and observe that
\[
\mu_\D^\o( \o\equiv -1\mbox{ on }p,\; \o\equiv-\epsilon \mbox{ on 
}\D\setminus p)>0
\]
whenever $\o(0,-n)=-1$.

(4) Uniform lower bounds for conditional probabilities were used 
twice, in the proofs of the shift lemma and the contour fluctuation 
lemma. In the proof of Lemma \ref{shift-lemma}, it is sufficient to replace 
the set $\{x,y\}$ by $\D(x)\cup\D(y)$, where 
$\D(x)=\{k-1,k,k+1\}\times\{n-1,n\}$ when $x=(k,n)$. This is because 
for $\o(k,n+1)=-1$ we have the estimate
\[
\mu_{\D(x)}^\o( \o(x)= -1,\; \o\equiv-\epsilon \mbox{ on 
}\D(x)\setminus \{x\})\geq \d \equiv 
\frac{\lambda\wedge1}{(1+\lambda)^6}\;.
\]

More care is needed in the proof of Lemma 
\ref{intersections} where we used a uniform estimate for the 
conditional probability of $B_{k,n}$ given $A_{k,n}$.
First, according to comment (5) on the proof of Theorem \ref{general} 
we have to specify the abscissas $a_n(\o),a_n(\o')$ by two 
parameters $k,k'\in\Z$ with $|k-k'|\leq 1$. Note, however, that 
the point $(a_n(\o),n)$ necessarily belongs to $\LL_{even}$  
because otherwise $\o=\epsilon$ at the adjacent points $(a_n(\o),n+1)$ 
and $(a_n(\o)+1,n+1)$; but this is excluded by the hard-core 
interaction. Therefore we have in fact $k=k'$, and we can consider 
the events $A_{k,n}$ as before. Next we 
redefine $\D_{k,n}$ as the set $\{k-1,\ldots,k+2\}\times \{n-1,n\}$,
and $B_{k,n}$ as the event that $\o(k,n)=\o(k+1,n)=1$ (as before),
$\o'(k,n)=\o'(k+1,n)=-1$ (in variation of the former definition),
and anything else occurs at the remaining sites of $\D_{k,n}$ (e.g., 
the vacuum configuration $-\epsilon$). We then have $d_n\geq2$ 
on $A_{k,n}\cap B_{k,n}$, and for $(\o,\o')\in A_{k,n}$ we find
\[
\mu^\o_{\D_{k,n}}\otimes \mu^{\o'}_{\D_{k,n}}(B_{k,n})\geq 
\frac{\lambda}{(1+\lambda)^8}\equiv\d
\]
as above. We can thus argue as before.

\end{document}